\documentclass[12pt]{amsart}

\usepackage{graphicx}
\usepackage{amssymb}
\usepackage{epstopdf}
\usepackage{amsfonts, amsmath, amssymb, amscd, latexsym,graphicx}
\newcommand{\Z}{{\bf Z}}
\newcommand{\Q}{{\bf Q}}
\newcommand{\N}{{\bf N}}

\newtheorem{theorem}{Theorem}
\newtheorem{corollary}[theorem]{Corollary}

\newtheorem{lemma}{Lemma}

\title{The word problem distinguishes counter languages}
\author{Sean Cleary}
\address{Department of Mathematics, The City College of New York \& The CUNY Graduate Center,
 New York, NY 10031, USA}
\email{cleary@sci.ccny.cuny.edu}
\author{Murray Elder}
\address{Department of Mathematics, Stevens Institute of Technology,
Hoboken, NJ 07030, USA }
\email{murrayelder@gmail.com}
\author{Gretchen Ostheimer}
\address{Department of Computer Science, Hofstra University, Hempstead NY 11549, USA}
\email{gretchen.ostheimer@hofstra.edu}
\thanks{The first author is grateful for the hospitality of the Centre de
Recerca Matem\`atica.}
\begin{document}
\maketitle

\begin{abstract}
Counter automata are more powerful versions of finite-state automata where addition and
subtraction operations are permitted on a set of $n$ integer registers, called counters.
We show that the word problem of $\Z^n$ is accepted by a nondeterministic $m$-counter automaton if 
and only if $m \geq n$.
\end{abstract}

\section{Introduction}

Connections between formal language theory and group theory have been considered 
by many authors. 
If $H$ is generated as a group by a finite set $X$, 
and if we let $X^\pm$ be the set $X$ together with formal inverses, 
one important language to consider is the {\em word problem}, which is
the set of words over $X^\pm$ representing the identity element of $H$. 
The formal language classification of the word problem of a group 
is independent of generating set in the sense that 
if ${\mathcal F}$ is a family of languages 
and if $X$ and $Y$ are two finite generating sets for a group $H$,
then the word problem of $H$ with respect to $X$
is in $\mathcal{F}$ if and only if the word problem of $H$ with respect to $Y$ is in ${\mathcal F}$
(see Gilman
\cite{Gilman:1996}).
Therefore we can refer to the word problem of a group 
rather than to the word problem of a particular generating set for a group.

It is natural then to ask about the extent to which the algebraic structure of a
group $H$ 
determines the formal language classification of the word problem of $H$.
In 1975 Anisimov and Seifert \cite{AnisimovSeifert:1975} proved that the word problem of $H$ is a regular
language if and only if $H$ is finite,
and in 1985 Muller and Schupp \cite{MullerSchupp:1983, MullerSchupp:1985} proved that the word problem of $H$ 
is a context-free language if and only if $H$ is virtually free. 
While the Anisimov and Seifert result can be proven easily from first principles, the Muller and Schupp 
result  relies heavily on a deep result of Stallings concerning one-ended groups \cite{stallings}.
In 1991 Herbst  \cite{Herbst:1991} used the Muller and Schupp result to show that the word problem of $H$ is a one-counter language if and only if $H$ is virtually cyclic. 
Notice that it follows from these results that 
if we restrict  our attention to languages which are word problems, 
nondeterministic automata which are either finite state, pushdown or one-counter are 
no more powerful than their deterministic counterparts. 

In formal language theory there are a variety of ways to generalize the ideas of finite-state, 
pushdown and one-counter automata. 
One such way is to consider $G$-automata, where $G$ is a group. 
Loosely, if $G$ is a group, a {\em
$G$-automaton} over a finite alphabet $X$ is an automaton in which
each edge is labeled by an ordered pair,
the first coordinate of which is an element of $G$
and the second coordinate of which is an element of $X^{\pm}$
or the empty word.
A word $w$ over $X^{\pm}$ is {\em
accepted by $A$} if there is a path from the initial state to a
final state for which the second coordinate reads the letters of  $w$ and the
product of the corresponding first coordinates is the identity element of $G$. 
If we take $G$ to be the trivial group, a $G$-automaton is simply a finite-state automaton,
and if we take $G = \Z$, a $G$-automaton is a one-counter automaton. 
For  $G = \Z^n$, a $G$-automaton is an $n$-counter automaton.
We show below  that the word problem of $\Z^n$ is accepted by a nondeterministic $m$-counter automaton if 
and only if $m \geq n$, so larger rank free abelian groups require more counters to accept
their word problems.  Thus, the natural heirarchy of counter languages coming from
the number of counters used does not collapse in the nondeterministic case of word problems of groups.

We note that sometimes counter automata are described as {\em blind counter automata}
(see Mitrana and Stiebe \cite{MitranaStiebe})
to emphasize the fact that the counters can not be examined until at an accept state.

A  pushdown-automaton is equivalent in power to 
a $G$-automaton where $G$ is free \cite{Kambites:overview}.
Kambites proved that for groups $G$ and $H$, 
$W(H)$ is accepted by a deterministic 
$G$ automaton if and only if $H$ has a finite index subgroup which embeds in $G$
\cite{Kambites:word}, so in the deterministic case, at least $n$ counters are
required to accept the word problem of  $\Z^n$.
Furthermore, he posed the following question: ``For what groups $G$ is it true that 
deterministic and non-deterministic $G$-automata accept the same word problems?"
Below, Theorem \ref{mainthm}  answers that question in the case that $G$ is abelian:
deterministic and non-deterministic $G$-automata accept the same word problems.
Our methods are elementary: we rely entirely on basic linear algebra. 

\section{Notation and definitions}

Let $G$ be a group. We define a {\em $G$-automaton} over $X$
to be a finite directed graph with a distinguished initial vertex,
some distinguished final vertices, and with edges labeled by $G
\times (X^{\pm}\bigcup\{\epsilon\})$
where $\epsilon$ is the empty word. 
We will refer to vertices as {\em states}.
By a {\em loop} we mean an edge that starts and ends at the same state,
and by a {\em circuit} we mean a path that does so. 

A $G$-automaton over $X$ is said to {\em accept} a word $w \in
X^{\pm*}$ if there is a path $p$ from the initial state to some
final state labeled $(1, w)$, where $1$ is the identity element of $G$.
In this case $p$ is called an {\em accepting path}.

If $n$ is a positive integer,
a $\Z^n$-automaton is called an $n$-{\em counter automaton}.
An $n$-counter language is one that is accepted by an $n$-counter automaton.
If $r$ is a regular expression over $X^{\pm}$, we let $L(r)$ denote the language
denoted by $r$. 

\section{General Preliminaries}

We will need to rely on two general results about languages accepted by $G$-automata.
The first establishes that having an $n$-counter word problem is a property of a group,
rather than of a particular generating set for the group.
The proof relies on basic properties of rational transductions as summarized 
by Gilman \cite{Gilman:1996} and Kambites \cite{Kambites:overview}, for example.
\begin{lemma}
\label{genSet}
If $G$ and $H$ are groups, and if $X$ and $Y$ are two finite generating sets for $H$, 
then the word problem for $H$ with respect to $X$ is accepted by a $G$-automaton 
if and only if the word problem for $H$ with respect to $Y$ is as well. 
\end{lemma}
\begin{proof}
Fix a group $G$, and 
let ${\mathcal F}$ be the set of languages which are accepted by some $G$-automaton.
Let $W_X$ be the word problem of $H$ with respect to $X$,
and let $W_Y$ be the word problem of $H$ with respect to $Y$.
Suppose that $W_X \in {\mathcal F}$. Then
$W_Y$  is a rational transduction of $W_X$
(see Proposition 2 in \cite{Kambites:overview}).
${\mathcal F}$ forms a family of languages. 
It follows that ${\mathcal F}$ is closed under rational transduction 
(see Theorem 6.2 in \cite{Gilman:1996}).
Therefore $W_Y \in {\mathcal F}$.
\end{proof}

The second general result establishes 
that the intersection of a regular language and an $n$-counter language is itself
$n$-counter. This is Lemma 3 in Elder \cite {Elder} and an immediate consequence of Theorem 4 in Kambites
 \cite{Kambites:overview}.
 Later we will need to refer to specific characteristics of an $n$-counter
 automaton that accepts such an intersection. 
For this reason we include the following lemma and proof: 
\begin{lemma}
\label{intersection}
Let $X$ be a finite set.
Let $L_1$ be regular language over $X$, and let $L_2$ be a language accepted
by a $\Z^n$-automaton over $X$. 
Then $L_1 \cap L_2$ is also accepted by a $\Z^n$-automaton over $X$.
\end{lemma}
\begin{proof}
Let $A_1$ be a finite-state automaton accepting $L_1$.
Let $A_2$ be a $\Z^n$-automaton accepting $L_2$.
We construct a $\Z^n$-automaton $B$ as follows.
The set of states of $B$ is $\Sigma_1 \times \Sigma_2$, where 
$\Sigma_i$ is the set of states of $A_i$. 
A state $(\sigma_1, \sigma_2)$ is final if and only if $\sigma_i$ is final 
in $A_i$ for $i=1,2$.
For $x \in X^{\pm}$ and $v \in \Z^n$, there is an edge $B$ from $(\sigma_1, \sigma_2)$ to $(\tau_1, \tau_2)$ labeled $(v, x)$ 
if and only if there is an edge in $A_1$ from $\sigma_1$ to $\tau_1$ labeled $x$
and there is an edge in $A_2$ from $\sigma_2$ to $\tau_2$ labeled $(v,x)$.
Furthermore, there is an edge  from $(\sigma_1, \sigma_2)$ to $(\tau_1, \tau_2)$
labeled $(v, \epsilon)$ if and only if one of three conditions holds:
\begin{itemize} 
\item there is an edge in $A_1$ from $\sigma_1$ to $\tau_1$ labeled $\epsilon$, and there is an edge in $A_2$ from $\sigma_2$ to $\tau_2$ labeled $(v, \epsilon)$, or
\item $v=0$, $\sigma_2 = \tau_2$ and there is an edge in $A_1$ from $\sigma_1$ to $\tau_1$ labeled $\epsilon$, or
\item $\sigma_1=\tau_1$ and there is an edge in $A_2$ from $\sigma_2$ to $\tau_2$ labeled $(v, \epsilon)$.
\end{itemize}
Words accepted by  $B$ are exactly those in  $L_1 \cap L_2$: a word $w$ in the intersection can follow a path labeled $(0,w)$ to states $(\sigma_1, \sigma_2)$ where $\sigma_i$ is a final state for $A_i$; similarly, any word accepted by $B$
can  lead to a state  $(\sigma_1, \sigma_2)$ which is a product of final states via a path
labeled $(0,w)$ and would thus be accepted by each of the $A_i$.
\end{proof}

\section{Main Result}

To show that we cannot accept the word problem of a free abelian group of
rank $n$ with a counter automaton with less than $n$ counters, we proceed
via a series of lemmas which allow us to consider automata of a preferred form
and to derive later a contradiction from a property somewhat analogous to the 
ranks of vectorspaces not being less than that of their subspaces.

We let $H = \Z^n$, and suppose that $x_1,x_2,\ldots,x_n$ is a basis for $H$ as 
a free abelian group. 
Let $X_1,X_2,\ldots,X_n$ be the formal inverses of the generators. 
Let $L = W(H) \cap L(x_1^*x_2^* \cdots x_n^*X_1^*X_2^* \cdots X_n^*)$.
If $j = (j_1,j_2,\ldots,j_n) \in \N^n$, 
let $w(j)$ denote the word $x_1^{j_1}x_2^{j_2 }\cdots x_n^{j_n}X_1^{j_1}X_2^{j_2} \cdots X_n^{j_n}$.

\begin{lemma}
\label{autStructure}
Let $H$ and $L$ be as above. 
Suppose $W(H)$ is $m$-counter. Then there is an $m$-counter automaton $A$ accepting
$L$ with the following structure:
\begin{itemize}
\item $A$ has a single final state $\sigma$. 
\item $A$ can be described as a collection of $2n$ subautomata \\
$A(x_1), A(x_2),\ldots, A(x_n), A(X_1), A(X_2), \ldots, A(X_n)$ satisfying the following 
criteria:
\begin{itemize}
\item the only edges between the subautomata are  labeled $(v, \epsilon)$ for some $v \in \Z^m$, and 
these edges go from $A(x_i)$ to $A(x_{i+1})$ for $i=1,2,\ldots,n-1$, from $A(X_i)$ to $A(X_{i+1})$ for $i=1,2,\ldots,n-1$, 
and from $A(x_n)$ to $A(X_1)$.
\item for all $a=x_i,X_i$, edges in $A(a)$ are labeled $(v,\epsilon)$ or  $(v,a)$ where $v \in \Z^n$.
\end{itemize}
\end{itemize}
\end{lemma}
\begin{proof}
Let $A_1$ be a finite-state automaton accepting the regular language $L(x_1^*x_2^* \cdots x_n^*X_1^*X_2^* \cdots X_n^*)$ of
the following specific form.
There are $n$ states $\sigma_{x_i}$, $n$ states $\sigma_{X_i}$, and two additional states $\alpha$, the initial state, and $\beta$, the only final state.
For $a=x_i,X_i$, the state $\sigma_a$ has a loop labeled $a$.
For $i=1,2,\ldots,n-1$, there are edges labeled $\epsilon$ from  $\sigma_{x_i}$ to $\sigma_{x_{i+1}}$ and from $\sigma_{X_i}$ to $\sigma_{X_{i+1}}$. In addition there are edges labeled $\epsilon$ from $\alpha$ to $\sigma_{x_1}$and from $\sigma_{X_n}$ to $\beta$.

Let $A_2$ be a $\Z^m$-automaton accepting $W(H)$.
We may assume without loss of generality that $A_2$ has a single final state. 
By Lemma \ref{intersection}, there exists a $\Z^m$-automaton $A$ accepting $L$. 
The automaton constructed in the proof of Lemma \ref{intersection}
has all of the desired properties. 
\end{proof}

We want to show that $m \geq n$ using linear algebra. The following lemma
will allow us to do so. 
${\N}$ denotes the set of positive integers. 
\begin{lemma}
\label{planes}
If $m < n$ then
$\N^n$ is not contained in the union of finitely many translates of  subspaces of $\Q^n$
each of which has dimension  at most $m$. 
\end{lemma}
\begin{proof}
Suppose that ${\bf N}^n$ is contained in the union of 
$Q_1, Q_2, \ldots, Q_r \subseteq \Q^n$,
where each $Q_i$ is a translate of an $m$-dimensional subspace of $\Q^n$.
Let $k = r+1$.
Let $B(k)$ be the set of all points $(x_1,x_2,\ldots, x_n)$  in ${\bf N}^n$ such that 
$x_i \leq k$ for $i=1,2,\ldots,n$. 
There are $k^n$ elements in $B(k)$. 
Let $B_i = B(k) \cap Q_i$. 
There are at most $k^m$ elements in $B_i$.
Therefore there at most $rk^m<k^{m+1} \leq k^n$ elements in $B(k)$.
We have reached a contradiction.  
\end{proof}

Let $p$ and $q$ be accepting paths in an $m$-counter automaton.
We will say that $p < q$ if $q$ can be obtained from $p$ by adding circuits. 
We will say that $p$ is {\em minimal} if it is minimal with respect to 
$<$.
\begin{lemma}
\label{paths}
Let $H$ and $L$ be as above. 
Suppose $W(H)$ is $m$-counter. Let $A$ be  an $m$-counter automaton $A$ accepting
$L$ with the structure posited in Lemma \ref{autStructure}.
There exist accepting paths $p$, $q_1, q_2, \ldots, q_n$ such that
\begin{itemize}
\item $p < q_i$ for $i=1,2,\ldots,n$;
\item  if $j \in \N^n$ such that $w(j)$ is the word accepted by $p$, 
and if $a_i \in \N^n$ such that $w(j+a_i)$ is the word accepted by $q_i$,
then $\{ a_1,a_2,\ldots,a_n \}$ is a set of linearly independent
vectors in $\N^n$.
\end{itemize}
\end{lemma}
\begin{proof}
Let $p$ be an accepting path which is minimal with respect to $<$,
and let $w(j)$ be the word that it accepts. 
Let $S_p$ be the semigroup spanned by all vectors of the form 
$j' - j$ such that there is a path $q$ accepting
$w(j')$ with $q > p$.
Consider the subspace $ V_p$ of ${\bf Q}^n$ spanned $S_p$.
Let $Q_p = j + V_p$.

There are finitely many accepting paths $p$ which are minimal with respect to $<$.
Suppose that none of these satisfies the criteria of the lemma.
Then each $V_p$ has dimension $n-1$ or smaller. 
But then ${\bf N}^n$ is contained in the union of finitely many
translates of subspaces  $\Q^n$ which are at most $(n-1)$-dimensional.
By Lemma \ref{planes} this is not possible. 
\end{proof}

\begin{theorem} 
\label{mainthm}
If the word problem of $\Z^n$ is an $m$-counter language, 
then $m \geq n$. 
\end{theorem}
\begin{proof}
Suppose that the word problem of $\Z^n$ with respect to some generating set is an 
$m$-counter language, with $m<n$. 
By Lemma \ref{genSet} we may assume that our generating set for $\Z^n$
is  a free basis $x_1,x_2,\ldots,x_n$.
Let $X_1,X_2,\ldots,X_n$ be formal inverses of the generators. 
Let $L = W(\Z^n) \cap L(x_1^*x_2^* \cdots x_n^*X_1^*X_2^* \cdots X_n^*)$.
By Lemma \ref{autStructure} there exists an $m$-counter automaton 
$A$ accepting $L$ with the specific structure
posited in that lemma.

We can take $p, j,q_i, a_i$ as in Lemma \ref{paths}.
Let $s_i$ be the $\Z^m$ contribution of the loops in $q_i$ that are not
in $p$ and which lie in $A(x_1) \cup A(x_2) \cup \cdots \cup A(x_n)$.
Let $S_i$ be the $\Z^m$ contribution of the loops in $q_i$ that are not
in $p$ and which lie in $A(X_1) \cup A(X_2) \cup \cdots \cup A(X_n)$.
Since $p$ and $q_i$ are both accepting, and since
$q_i$ is built up from $p$ in the specific way that it is, 
$s_i + S_i = 0$. 

Since any set of $n$ vectors in $\Z^m$ is linearly dependent,  then
there exist $\alpha_1, \alpha_2, \ldots, \alpha_n \in \Z$ not all zero
such that $\alpha_1 s_1 + \alpha_2 s_2 + \cdots + \alpha_n s_n = 0$.
We construct an accepting path $r$ as follows.
We start with $p$. 
If $\alpha_i$ is strictly positive, 
consider those loops of $q_i$ that are not part of $p$ 
but that do lie in $A(x_1) \cup A(x_2) \cup \cdots \cup A(x_n)$;
add $\alpha_i$ times as many traversals of  these loops.
The $\Z^m$ contribution of these loops is 
$\alpha_i s_i$.
If $\alpha_i$ is strictly negative,
do the same thing but this time 
consider those loops of $q_i$ that are not part of $p$ 
but that do lie in $A(X_1) \cup A(X_2) \cup \cdots \cup A(X_n)$,
and add $- \alpha_i$ times as many traversals of  these loops. 
The $\Z^m$ contribution of these loops is 
$(- \alpha_i )S_i = \alpha_i s_i$. 
The path $r$ is accepting since 
$\alpha_1 s_1 + \alpha_2 s_2 + \cdots + \alpha_n s_n = 0$.

We now reach a contradiction by showing that the word accepted by $r$ does
not represent the identity and thus is not in $W(\Z^n)$.
Consider the case, for example, when $\alpha_1 < 0$ and 
$\alpha_2, \alpha_3, \ldots, \alpha_n \geq 0$.
Let $u = j + \alpha_2 a_2 + \cdots + \alpha_n a_n$, and 
let $v = j + (-\alpha_1) a_1$.
Then the word $w$ accepted by $r$ is of the form
$x_1^{u_1}x_2^{u_2} \cdots x_n^{u_n}X_1^{v_1}X_2^{v_2} \cdots X_n^{v_n}$, so 
$w$ is in the word problem only if $u = v$. 
This is the case if and only if
\begin{eqnarray*}
\alpha_2 a_2 + \alpha_3 a_3 + \cdots + \alpha_n a_n& =&-\alpha_1a_1
\end{eqnarray*}
This is impossible since $\{a_1, a_2, \ldots, a_n \}$ is linearly independent. 
All other cases reach a similar contradiction. 
\end{proof}
From the classification of finitely-generated abelain groups, we get the immediate corollary, analgous to Kambites Theorem 1 \cite{Kambites:word} for
the group case but in the nondeterministic case:
\begin{corollary}
\label{generalabelian}
The word problem of finitely-generated abelian group $H$ is recognized by
a nondeterministic $G$-automaton if and only if $H$ has a finite-index subgroup isomorphic
to a subgroup of $G$.
\end{corollary}

\section{Acknowledgments}

We would like to thank Mark Kambites for the relevant background material from 
semigroup theory and Bob Gilman for suggesting the proof of Lemma \ref{planes}.

\bibliographystyle{plain}

 \end{document}